\documentclass[12pt]{article}
\usepackage{amsmath,amssymb,amsfonts,amsthm}
\usepackage{eucal} 
\usepackage{graphicx}
\usepackage[small,nohug,heads=littlevee]{diagrams}
\diagramstyle[labelstyle=\scriptstyle]

\newcounter{ecount}
\newcommand{\C}{\mathbb{C}}

\newcommand{\F}{\mathbb{F}_3}

\newcommand{\Fp}{\mathbb{F}_p}
\newcommand{\Fpd}{\Fp^d}

\newcommand{\bm}{\mathbf{m}}

\newcommand{\g}{\mathfrak g}

\newcommand{\Ug}{\mathcal{U}(\g)}
\newcommand{\zee}{{\mathfrak z}}

\DeclareMathOperator{\End}{End}

\DeclareMathOperator{\Tr}{Tr}

\DeclareMathOperator{\diag}{diag}
\DeclareMathOperator{\Aut}{Aut}

\DeclareMathOperator{\eval}{eval}

\DeclareMathOperator{\rank}{rk_{\Fp}}

\newtheorem{Theorem}{Theorem}

\author{Erik Carlsson}

\title{Symmetric Functions and Caps}

\begin{document}

\maketitle

\begin{abstract}
Given a finite subset $S \subset \Fpd$, let $a(S)$ be the number of distinct $r$-tuples
$\alpha_1,...,\alpha_r \in S$ such that $\alpha_1+\cdots+\alpha_r = 0$. We consider
the ``moments'' $F(m,n) = \sum_{|S|=n} a(S)^m$.
Specifically, we present an explicit formula for $F(m,n)$ as a product of two matrices,
ultimately yielding a polynomial in $q=p^d$. The first matrix is independent of $n$ while the second
makes no mention of finite fields. However, the complexity of calculating each grows with $m$. 
The main tools here are the Schur-Weyl duality theorem, and some elementary
properties of symmetric functions. This problem is closely to the study of maximal caps.

\end{abstract}

\section{Introduction}

Given a subset $S \subset \Fpd$, denote by $a(S)$ the number of distinct $r$-tuples such that 
$\sum_j \alpha_j = 0$ in $\Fpd$. An important ``zero-sum'' problem \cite{GG} is determining how
large $|S|$ can be subject to the constraint $a(S) = 0$. In this paper, we consider the functions
\[F(m,n) = \sum_{|S|=n} a(S)^m,\]
where $S \subset \Fpd$, $p$ is a prime. $F(m,n)/F(0,n)$ is the $m$th 
moment of the push-forward of the uniform distribution through $a$. Since this distribution is
compactly supported, the values $F(m,n)$ for $0 \leq m \leq \max_{|S|=n} a(S)$ determine the
entire distribution.

The strategy of this paper is to describe $F(m,n)$ as the trace
of an operator that does not depend on $n$. Sections \ref{decomp} and \ref{algelt}
describe the operator, and the decomposition mentioned in the abstract.
Section \ref{thmsect} is the main theorem about the decomposition, while
sections \ref{g} and \ref{h} give explicit formulas for the matrix elements. Finally, section
\ref{f} covers the example $m=2$, $n=10$, $p=3$, $d$ arbitrary.

\section{The Moment Function as a Trace}

\label{decomp}

Let $V = \C \cdot \Fpd$, the $\C$-vector space with basis
$v_\alpha, \alpha \in \Fpd$, $\g = \End(V)$, the Lie algebra of $GL(V)$, and
\[T(\g) = \bigoplus_{k \geq 0} T_k = \bigoplus_{k \geq 0} \g^{\otimes k},\]
the tensor algebra of $\g$. For any representation
$\rho : GL(V)\rightarrow GL(\rho(V))$ we have a map
\[\varphi=\varphi_\rho : T(\g) \rightarrow \Ug \rightarrow \End(\rho(V))\]
given by
\[\varphi(A_1 \otimes \cdots \otimes A_k) = \rho'(A_1) \circ \cdots \circ \rho'(A_k),\]
where $\rho'$ is the derivative of $\rho$.

When $\rho = \wedge^n$, $\rho(V)$ is the vector space with basis $v_S$, $\#S = n$.
Our aim is to represent the operator
\[v_S \mapsto a(S) v_S\]
as $\varphi(B)$ for some $B\in T(\g)$. Then $F(m,n)$ takes the form
\[F(m,n) = h(B^{\otimes m}),\quad h=\Tr\circ \varphi_{\wedge^n}.\]
It turns out that one can choose $B$ so that $\deg B \leq p$, and so that $B$ is generated by
simultaneously commuting elements of $\g$. One can therefore take $F(m,n) = h(B^m)$, where
$B^m = \Phi\left(B^{\otimes m}\right)$, and $\Phi=\Phi_k$ is the symmetrization map
\[\Phi (A_1\otimes \cdots \otimes A_{k}) = \frac{1}{k!}\sum_{\sigma \in S_{k}} 
\sigma \cdot A_{\sigma(1)}\otimes \cdots \otimes A_{\sigma(k)}.\]

We capitalize on this element using a factorization of the composed map $f = h \circ \Phi$:
$f$ is a $GL(V) \times S_{k}$-invariant map to $\C$, so factors through
\[f: \g^{\otimes k} \xrightarrow{g} \left(\g^{\otimes k}\right)_{trivial} \xrightarrow{h} \C,\]
where the middle term is the trivial subspace with respect to the $GL(V) \times S_{k}$ action. 

Our theorem is a calculation of the matrix elements of $g$, $h$ in a canonical basis of
$\left(\g^{\otimes k}\right)_{trivial}$. The use of this decomposition is
can be seen from a complexity point of view: the matrix elements of $g$
depend on $\Fpd$, and $m$, but do not contain an $n$ term. The complexity of computing $h$
depends on $m$ and $n$, but not on $B$.

\section{The Tensor Algebra Element}

\label{algelt}

To find $B$, we use a fact used by Bierbrauer and Edel in \cite{BE}: if $\alpha \in \Fpd$,
then
\[\frac{1}{q}\sum_{\beta \in \Fpd} \zeta^{\alpha\cdot \beta} = \delta_{\alpha,0},
\quad \zeta = e^{2\pi i/p}.\]
If $S = \left\{\alpha_1,...,\alpha_n\right\}$, then
\begin{equation}
\label{be}
a(S) = \frac{1}{q}\sum_\beta e_r(\zeta^{\alpha_1\cdot\beta},...,\zeta^{\alpha_n\cdot\beta})
\end{equation}
where $e_r$ is the $r$th elementary symmetric polynomial. $e_r$ has a useful expression in the power-sum basis
\begin{equation}
\label{en2pmu}
e_r = \sum_{\rho} \frac{\langle e_r,p_\rho\rangle}{\zee(\rho)} p_\rho,
\end{equation}
where $\rho \in \Lambda(r)$, the set of partitions of $r$, $p_\rho = \prod_{j=1}^{\ell=\ell(\rho)} p_{\rho_j}$, $p_k = \sum_i x_i^k$, and
\[\langle p_\rho,p_{\rho'}\rangle = \delta_{\rho,\rho'}\zee(\rho),\quad
\zee(\rho) = |\Aut(\rho)|\prod_j \rho_j,\]
the standard inner product on symmetric polynomials.

We can now describe $B$. Let
\[A_\beta:V\rightarrow V,\quad A_\beta \cdot v_\alpha = \zeta^{\alpha\cdot\beta} v_\alpha,\]
and
\begin{equation}
\label{mdecomp}
B = \frac{1}{q}\sum_{\beta,\rho} \frac{\langle e_r,p_\rho\rangle}{\zee(\rho)} 
A_{\rho_1\beta} \cdots A_{\rho_\ell \beta}
\end{equation}
in $T_{k\leq p}$. The product here is the symmetric product, i.e. the symmetrization of the tensor product.
By (\ref{be}) and (\ref{en2pmu}), $B$ has
the desired property $f(B^m) = F(m,n)$.

\section{Schur-Weyl Duality and The Theorem}

\label{thmsect}

The canonical basis of $\left(\g^{\otimes k}\right)_{trivial}$ is better said
under the isomorphism
\begin{equation}
\label{Vident}
\End(V)^{\otimes k} \cong V^{\otimes k} \otimes V^{*\otimes k} \cong \End(V^{\otimes k}).
\end{equation}
The decomposition is now Schur's lemma and Schur-Weyl duality:
\[\left(\g^{\otimes k}\right)_{trivial} \cong \bigoplus_{|\mu| = k} \C \cdot P_{\mu}.\]
Here $\mu$ is a partition of $k$, and $P_{\mu} \in \End(V^{\otimes k})$ is the projection onto 
$\mathbb{S}_\mu(V)\otimes S_\mu$ in the Schur-Weyl decomposition
\[V^{\otimes k} \cong \bigoplus_{\mu} \mathbb{S}_\mu(V)\otimes S_\mu\]
into irreducible representations of $GL(V)\times S_k$. $\mathbb{S}_\mu$ is Schur functor, and $S_\mu$ is the irreducible Specht module  \cite{Sa} associated to $\mu$.
The fact that the multiplicity of each irreducible is a most one means that $P_\mu$ is a canonical basis.
We calculate the matrix elements of $g$, $h$ in this basis.

\begin{Theorem}
The moment function is given by
%
%
%
\[F(m,n) = \sum_{|\mu| = k} 
\frac{G(m,\mu)\cdot H(\mu,n)}{\dim \mathbb{S}_\mu(V) \cdot \dim S_{\mu}},\]
where
\[G(m,\mu) = \Tr_{V^{\otimes k}} B^m \circ P_{\mu},\quad
H(\mu,n) = h(P_{\mu}).\]
Furthermore, combinatorial formulas for $G$, $H$, $F$ are given in (\ref{G}), (\ref{H}) and (\ref{F}),
and the answer is a polynomial of degree $\leq n-1$ in $q$.

\end{Theorem}

\begin{proof}
This just expresses $f(B^m) = h\cdot g(B^m)$ in the basis $P_\mu$. The explicit formulas
for $G$, $H$, and $F$ are given in sections \ref{g}, \ref{h}, and \ref{f}. The fact that
the answer is a polynomial in $q$ of degree $n-1$ follows from (\ref{F}), and the fact that
$m_\nu(\lambda) = 0$ for $\ell(\lambda) < \ell(\nu)$.
The bound for the degree also follows from the simple fact that
\[\lim_{q\rightarrow \infty} \frac{F(m,n)}{F(0,n)} = 0,\]
so that $\deg(F(m,n)) < \deg(F(0,n)) = n$.

\end{proof}

\section{Decomposition of the Algebra Element}
\label{g}

This section calculates $G(m,\mu)$, the coordinate of $g(B^m)$.
First, we expand $B^m$:
\begin{equation}
\label{Bbm}
B^m = \sum_{|\bm|=m} \prod_\rho \left(\frac{\langle e_r,p_\rho\rangle}{\zee(\rho)}\right)^{\bm(\rho)}
\cdot \frac{m!}{\prod_\rho \bm(\rho)!} \cdot \Phi\left(B_\bm\right),
\end{equation}
Where the sum is over functions $\bm: \Lambda(r) \rightarrow \mathbb{N}$ with $|\bm| = \sum_\rho \bm(\rho) = m$, 
$k = \sum_\rho \bm(\rho)\ell(\rho)$, and
\[B_\bm = \frac{1}{q^m}\bigotimes_{\rho} 
\left(\sum_\beta A_{\rho_1\beta} \otimes \cdots \otimes A_{\rho_\ell\beta}\right)^{\otimes \bm(\rho)} \in T_k,\]
with the tensor taken in the lexicographic order on $\Lambda(r)$.
Let $G(\bm,\mu)$ be the coordinate of each component $g(B_\bm)$,
\[G(\bm,\mu) =  \Tr_{V^{\otimes k}} B_\bm \circ P_\mu.\]

Given $|\mu| = |\nu| = k$, let $\chi^\mu_\nu$ denote the character of the irreducible 
representation of $S_{k}$ corresponding to $\mu$ on $\sigma \in S_{k}$ of cycle-type $\nu$. 
By Schur-Weyl duality, this is the same as $\langle s_\mu,p_\nu\rangle$, where $s_\mu$ is the Schur polynomial.
Using \cite{Se}, theorem 8(ii),
\[G(\bm,\mu) = \sum_\nu  \frac{n_\mu\chi^\mu_\nu}{k!} G_0(\bm,\nu),\quad
G_0(\bm,\nu) = \sum_{\sigma \sim \nu} 
\Tr \bigotimes_\rho B_\bm \circ \sigma,\quad n_\mu = \dim S_\mu.\]

Now let $\pi$ denote a function
\[\left\{1,...,k\right\}\rightarrow \left\{1,...,\ell\right\},\quad \ell = \ell(\nu),\]
and $[\pi]$ the corresponding partition of $\left\{1,...,k\right\}$ into the
level sets of $\pi$. Associated to each $\sigma$ we have a partition $[\pi](\sigma)$ 
which groups numbers if they are in the same cycle in $\sigma$. Also set $[\pi_\bm]$
denote the set partition corresponding to $\bm$:
\[\left\{1,...,k\right\} = \bigsqcup_\rho \bigsqcup_{1\leq a \leq \bm(\rho)} C_{\rho,a},\] 
with $C_{\rho,a}$ a block of size $\ell(\rho)$, and where the disjoint union is with respect
to the lexicographic ordering of $\Lambda(r)$. This associates to each $1\leq c\leq k$ a $3$-tuple
\begin{equation}
\label{k2rhoji}
c \leftrightarrow (\rho,a,b),
\end{equation}
$b$ being the index within $C_{\rho,a}$.

The summand in $G_0$ only depends on $[\pi](\sigma)$:
\[G_0(\bm,\nu) = \sum_{[\pi] \sim \nu} \#\left\{\sigma, [\pi](\sigma) = [\pi]\right\} \]
\[\sum_{\alpha_1,...,\alpha_{k}}\delta(\alpha_c \sim \alpha_{c'}\mbox{ if $c \sim c'$ in $[\pi]$})
\delta(\sum_{c \in C_{\rho,a}} \rho_b \alpha_c \cong 0, \mbox{for each $\rho,a$}).\]
The sum over $\alpha_j$ is just the dimension of the kernel over $\Fp$ of the $m$ by $\ell$ matrix 
$X = X(\bm,\pi)$, where
\[X_{i,j} = \sum_{c \in \pi_\bm^{-1}(i) \cap \pi^{-1}(j)} \rho_b,\]
with $\rho,a,b$ associated to $c$ as in (\ref{k2rhoji}).

The summand depends only on the equivalence class $[X]_{\bm,\ell}$, where
\[[X]_\bm = \prod_{\rho} S_{\bm(\rho)} \cdot X,\quad 
[X]_\ell = X \cdot S_{\ell},\quad [X]_{\bm,\ell} = [[X]_\bm]_\ell.\]
Now $G_0$ becomes
\[G_0(\bm,\nu) = \sum_{[X]_{\bm,\ell}}G_1(\bm,[X]_{\bm,\ell},\nu) q^{\dim \ker X\ (\mbox{\tiny mod $p$})}\]
\[G_1(\bm,[X]_{\bm,\ell},\nu) = \#\{\sigma, \sigma \sim \nu, [X]_{\bm,\ell}(\sigma) = [X]_{\bm,\ell}\}.\]
It is helpful to think of $[X]_{\bm,\ell}$ as an equivalence class of the bipartite graph with adjacency matrix $X$.

Since every map is onto with fibers of the same size in
\begin{diagram}
&& \pi & \rMapsto^{g_2} & [X]_\ell \\
&& \dMapsto^{g_3} & & \dMapsto^{g_4}\\
\sigma \sim \nu & \rMapsto^{g_1}& [\pi] & \rMapsto^{g_5} & [X]_{\bm,\ell}
\end{diagram}
we get
\[G_1(\bm,[X],\nu) = \frac{|g_1^{-1}([\pi])|\cdot |g_2^{-1}([X]_\ell)|\cdot |g_4^{-1}([X]_{\bm,\ell})|}{|g_3^{-1}([\pi])|} =\]
\[\frac{\prod_j (\nu_j-1)! \#[X]_{\bm,\ell}}
{\ell!} |g_2^{-1}([X]_\ell)|,\]
\[|g_2^{-1}([X]_{\ell})| = \sum_{X \in [X]_\ell} \#\left\{\pi\sim \nu| \pi \mapsto X\right\} =
\frac{\ell!}{|\Aut(\nu)|}[m_\nu] \Phi_\ell\left(J(X)\right),\]
where $[\mu_{\nu}]$ is the coefficient of $m_{\nu}$,
\[J(X)(x_1,x_2,...) = \prod_i \sum_{\pi_0} 
\delta\left(X_{ij}=\sum_{b \in \pi_0^{-1}(j)} {\rho_b},\mbox{ each $j$}\right) \prod_b x_{\pi_0(b)},\]
$\rho$ is associated to $1\leq i \leq m$ under $\pi_\bm$,
and the sum is over $\pi_0 :\left\{1,...,\ell(\rho)\right\} \rightarrow \left\{1,...,\ell\right\}$. The final formula is
%
%
\[G(\bm,\mu) = \sum_{[X]_{\bm,\ell}} \#[X]_{\bm,\ell}\cdot q^{\dim \ker X\ (\mbox{\tiny mod $p$})}\cdot\]
\begin{equation}
\label{G}
%
%
\sum_\nu \left(\frac{n_\mu \chi^\mu_\nu}{k!\zee(\nu)} \prod_j \nu_j!\right) \cdot [m_\nu] \Phi_\ell\left(J(X)\right).
\end{equation}
\section{Trace of the Schur Projections}
\label{h}
Finally, we calculate the matrix element
\[H(\mu,n) = \Tr \varphi_{\wedge^n}\left(P_{\mu}\right) =\]
\begin{equation}
\label{H02H}
\sum_{\nu}\frac{n_\mu\chi^\mu_\nu}{k!} H_0(\nu,n),\quad H_0(\nu,n) = \sum_{\sigma\sim\nu} \Tr \varphi_{\wedge^n}(\sigma),
\end{equation}
applying $\varphi$ using the identification (\ref{Vident}).

Proceeding directly leads to difficult combinatorics, but the computation can be much simplified using
the identification of symmetric polynomials with (virtual) representations of
$GL(V)$. Specifically, one identifies
\[\wedge^n(V) \leftrightarrow e_n,\]
$e_n(x_1,...,x_q)$ being the character of $\wedge^n (\diag(x_1,...,x_q))$. By (\ref{en2pmu}),
\begin{equation}
\label{H12H0}
H_0(\nu,n) = \sum_\lambda \frac{\langle e_n,p_\lambda\rangle}{\zee(\lambda)}H_1(\nu,\lambda),\quad
H_1(\nu,\lambda) = \sum_{\sigma \sim \nu} \Tr \varphi_{\rho_{\lambda}}(\sigma),
\end{equation}
where $\rho_{\lambda}(V)$ is the vector space with basis 
\[v_I = v_{i_1}\otimes \cdots \otimes v_{i_\ell}, \quad \ell = \ell(\lambda),\]
$v_i$ are basis vectors of $V$, and
\[\rho_{\lambda}(x) v_I = \left(x^{\lambda_1} v_{i_1}\right) \otimes \cdots \otimes \left(x^{\lambda_\ell} v_{i_\ell}\right).\]
$\rho_\lambda$ is not representation, but 
\begin{equation}
\label{phirhod}
\rho_\lambda'(\xi) v_I = \sum_k v_{i_1} \otimes \cdots \otimes \left(\lambda_k  \xi\cdot v_{i_k}\right) \otimes \cdots \otimes v_{i_\ell}
\end{equation}
is a linear map $\End(V) \rightarrow \End(\rho_\lambda(V))$ so $\varphi_{\rho_\lambda}$ is well-defined.
$H_1$ turns out to be the better-behaved expression.

Continuing,
\[H_1(\nu,\lambda) = \sum_{\sigma \sim \nu} \Tr \sum_J \varphi_{\rho_\lambda}
(E_{\sigma(j_1),j_1} \otimes \cdots \otimes E_{\sigma(j_{k}),j_{k}})=\]
\[\sum_{\sigma \sim \nu} \sum_{J,K} \langle\rho_\lambda'
(E_{\sigma(j_1),j_1}) \circ \cdots \circ \rho_\lambda'(E_{\sigma(j_{k}),j_{k}})v_K,v_K\rangle, \]
\[\langle v_I,v_J \rangle = \delta_{I,J}.\]

We calculate each summand,
\[\langle\rho_\lambda'(E_{i_1,j_1}) \circ \cdots \circ \rho_\lambda'(E_{i_{k},j_{k}})v_K,v_K\rangle =\]
\[\sum_\pi \prod_{1\leq b\leq \ell} \lambda_b^c\langle E_{i_{a_1},j_{a_1}} \circ \cdots
\circ E_{i_{a_c},j_{a_c}} v_{k_b},v_{k_b}\rangle =\]
\begin{equation}
\label{sortedcoefficient}
\sum_\pi \prod_{1\leq b\leq \ell} \lambda_b^c\delta(k_b-i_{a_1})\delta(j_{a_1}-i_{a_2})\cdots\delta(j_{a_c}-k_b),
\quad \delta(k) = \delta_{k,0}.
\end{equation}
The sum is over $\pi:\{1,...,k\}\rightarrow \{1,...,\ell\}$ which keeps track of where each $E_{i_a,j_a}$
is inserted as in (\ref{phirhod}). Also, we have set $\{a_1,...,a_c\} = \pi^{-1}(b)$, with the ordering $a_1 < \cdots < a_c$.

Associated to each $\pi$ there is a unique permutation $\sigma'(\pi) \in S_{k}$ which is ``sorted,'' and such that
$[\pi](\sigma') = [\pi]$. Sorted here means that each cycle of $\sigma'$ can be written as $(a_1,...,a_c)$ with
$a_1 < \cdots < a_c$. The summand depends only on $\sigma'(\pi)$. 
Summing over $J$, $K$ and setting $\sigma' \sim \nu'$, $\sigma^{-1}\cdot {\sigma'} \sim \nu''$ gives
\[H_1(\nu,\lambda) = \sum_{\sigma \sim \nu,\sigma' \mbox{\tiny -sorted}} q^{\ell(\lambda)-\ell(\nu')}
\left(\sum_{\pi,[\pi] = [\pi](\sigma')}\prod_b \lambda_b^c\right)
\sum_{j_1,...,j_{k}} \prod_k\delta(j_{\sigma'(k)}-j_{\sigma(k)}).\]
Since the summand is invariant under simultaneously conjugating
\[(\sigma,\sigma') \mapsto (\tau\sigma\tau^{-1},\tau \sigma' \tau^{-1}),\]
we may exchange the ``sorted'' condition for a factor of $\prod_j (\nu'_j-1)!$:
\[\sum_{\sigma \sim \nu,\sigma'}\frac{|\Aut(\nu')|}{\prod_j (\nu'_j-1)!}
q^{\ell(\lambda)-\ell(\nu')+\ell(\nu'')} m_{\nu'}(\lambda) =\]
\[\sum_{\nu',\nu''} \frac{\zee(\nu')}{\prod_j \nu'_j!}
q^{\ell(\lambda)-\ell(\nu')+\ell(\nu'')}H_2(\nu,\nu',\nu'')(\lambda),\]
\[H_2(\nu,\nu',\nu'') = \#\left\{\sigma\sim\nu,\sigma'\sim\nu'|\sigma\cdot\sigma'\sim\nu''\right\},\]
and $m_{\nu'}$ is the monomial symmetric function.

This expression can be simplified:
first, $H_2$ is the same as the coefficient of $\sum_{\sigma'' \sim \nu''} \sigma''$ in 
\[\left(\sum_{\sigma \sim \nu} \sigma\right)\cdot\left(\sum_{\sigma'\sim \nu'} \sigma'\right)\]
in the center of the group-ring $\C[S_{k}]$. Passing to the idempotent basis, one easily deduces that
\begin{equation}
\label{H2simple}
H_2(\nu,\nu',\nu'') = \frac{(k!)^2}{\zee(\nu)\zee(\nu'')\zee(\nu'')}
\sum_\rho \frac{\chi_{\nu}^\rho\chi_{\nu'}^\rho\chi_{\nu''}^\rho}{n_\rho}.
\end{equation}
Summing over $\nu''$ gives
\begin{equation}
\label{H1simple}
H_1(\nu,\lambda) = \frac{(k!)^2}{\zee(\nu)}\sum_{\nu'}\frac{1}{\prod_j \nu'_j!} 
\sum_\rho \dim \mathbb{S}_\rho\cdot \frac{\chi_{\nu}^\rho\chi_{\nu'}^\rho}{n_\rho}m_{\nu'}(\lambda),
\end{equation}
where we have used (\ref{H2simple}), and
\[\sum_{\nu''} \frac{q^{\ell(\nu'')}}{\zee(\nu'')}\chi^\rho_{\nu''} = \eval(s_\rho,p_k \mapsto q) = 
\dim \mathbb{S}_\rho(V).\]
Finally, summing over $\nu$ gives
\begin{equation}
\label{H}
H(\mu,n) = 
%
k!\dim \mathbb{S}_\mu(V)\sum_{\nu',\lambda}q^{\ell(\lambda)-\ell(\nu')}
\frac{\chi^{(1^n)}_{\lambda} \chi_{\nu'}^{\mu}}
{\zee(\lambda)\prod_j \nu'_j! }
m_{\nu'}(\lambda),
\end{equation}
using (\ref{H1simple}), (\ref{H02H}), (\ref{H12H0}), and orthogonality of characters.

\section{An Example}
\label{f}

Combining equations (\ref{Bbm}), (\ref{G}), (\ref{H}), and the formula
\[\sum_\mu \chi^\mu_\nu\chi_{\nu'}^\mu = \delta_{\nu,\nu'} \zee(\nu),\]
we arrive at a formula for $F$:
\[F(m,n) = \sum_\lambda \left(\frac{q^{\ell(\lambda)}\chi^{(1^n)}_\lambda}
{\zee(\lambda)}\right) {\mathbf{F}}(m)(\lambda),\]
where
\[{\mathbf{F}}(m) = \sum_{|\bm| = m} \frac{m!}{\prod_\rho \bm(\rho)!}
\prod_\rho \left(\frac{\chi^{(1^r)}_\rho}{\zee(\rho)}\right)^{\bm(\rho)} \mathbf{F}(\bm),\]
\begin{equation}
\label{F}
\mathbf{F}(\bm)(x_1,x_2,...) = \sum_{[X]_{\bm,\ell}} q^{-\rank X}\#[X]_{\bm,\ell}\cdot \Phi_\ell \left(J(X)\right).
\end{equation}
$\mathbf{F(\bm)}$ is a symmetric polynomial of degree $k$, independent of $n$.

As an example, we present the second moment for $n=10$ points in $\F^d$. 
As the formula shows, this requires
a sum over partitions of size $k\leq 6$, equivalence classes of bipartite graphs with $\leq 6$ edges, and over 
the character table of $S_{10}$. The answer is
\[\frac{F(2,10)}{F(0,10)} = \frac{120(q^2+89q-540)}{(q-5)(q-4)(q-2)}.\]
Notice that $F(2,10)$ vanishes at $q=1,3,9$, since there are no subsets of size $10$
when $d=0,1,2$.


\begin{thebibliography}{99}

\bibitem{BE}
Bierbrauer, J\"{u}rgen; Edel, Yves, \emph{Bounds on affine caps}. J. Combin. Des. 10 (2002), no. 2, 111--115. (Reviewer: Tamás Sz\"{o}nyi) 51E22 (05B25).

\bibitem{DM}
D. Maclagan and B. Davis, \emph{The Card Game SET! The Mathematical Intelligencer}. Volume 25, Number 3, 2003.

\bibitem{FH}
Fulton,William and Harris, \emph{Joe, Representation Theory, A First Course},
Springer, 1991

\bibitem{GG}
W. Gao and A. Geroldinger, \emph{Zero-sum problems in finite abelian groups: a survey}.
Expo. Math., 24:337–369, 2006

\bibitem{Mac}
I.~Macdonald, 
\emph{Symmetric functions and Hall polynomials}, 
The Clarendon Press, Oxford University Press, New York, 1995.

\bibitem{Sa}
B.E. Sagan, \emph{The Symmetric Group: Representations, Combinatorial Algorithms,
and Symmetric Functions, 2nd Ed}, Springer-Verlag, New York, 2001.

\bibitem{Se}
J.-P. Serre, \emph{Linear Representation of Finite Groups}, Springer-Verlag, New York, 1977.

\end{thebibliography}
\end{document}